\documentclass[a4paper]{article}
\usepackage{latexsym}
\usepackage{amssymb,amsmath}
\usepackage{color}

\DeclareMathOperator{\rk}{rk}
\DeclareMathOperator{\tr}{tr}

\newtheorem{Theorem} {Theorem} [section]
\newtheorem{Proposition} [Theorem] {Proposition}
\newtheorem{Corollary} [Theorem] {Corollary}
\newcommand{\Proof}{ \noindent{\bf Proof.}\quad }
\newcommand{\qed}{\hfill$\Box$\medskip}
\newcommand{\adj}{\sim}
\newcommand{\maysplit}{\discretionary{}{}{}}
\title{Distance-regular graphs where the distance-$d$ graph
has fewer distinct eigenvalues}
\author{A. E. Brouwer \& M. A. Fiol$^{*1}$ \vspace{8pt} \\
 \vspace{-3pt}
 \small $^*$Universitat Polit\`ecnica de Catalunya, BarcelonaTech \\
 \vspace{-3pt}
 \small Dept. de Matem\`atica Aplicada IV, Barcelona, Catalonia \\
 \small (e-mails: {\tt aeb@cwi.nl, fiol@ma4.upc.edu})
}

\date{August 5, 2014}
\begin{document}
\maketitle
\footnotetext[1]{\hbadness=1050
Research supported by the
{\em Ministerio de Ciencia e Innovaci\'on}, Spain, and the
{\em European Regional Development Fund} under project MTM2011-28800-C02-01,
and the {\em Catalan Research Council} under project 2009SGR1387.
}

\begin{abstract}
Let the Kneser graph $K$ of a distance-regular graph $\Gamma$
be the graph on the same vertex set as $\Gamma$, where two vertices
are adjacent when they have maximal distance in $\Gamma$.
We study the situation where the Bose-Mesner algebra of $\Gamma$
is not generated by the adjacency matrix of $K$.
In particular, we obtain strong results in the so-called `half antipodal' case.
\end{abstract}
\medskip

\medskip


\bigskip
Let $\Gamma$ be a distance-regular graph of diameter $d$ on $n$ vertices.
Let $\Gamma_i$ be the graph with the same vertex set as $\Gamma$
where two vertices are adjacent when they have distance $i$ in $\Gamma$.
Let $A$ be the adjacency matrix of $\Gamma$, and $A_i$ that of $\Gamma_i$.
We are interested in the situation where $A_d$ has fewer distinct
eigenvalues~than~$A$.
In this situation the matrix $A_d$ generates a proper subalgebra of the
Bose-Mesner algebra of $\Gamma$, a situation reminiscent of imprimitivity.
We survey the known examples, derive parameter conditions,
and obtain strong results in what we called the `half antipodal' case.
Unexplained notation is as in \cite{BCN}.

\medskip
%
%
The vertex set $X$ of $\Gamma$ carries an association scheme
with $d$ classes, where the $i$-th relation is that of having
graph distance $i$ ($0 \le i \le d$).
All elements of the Bose-Mesner algebra ${\cal A}$
of this scheme are polynomials of degree at most $d$ in the matrix $A$.
In particular, $A_i$ is a polynomial in $A$ of degree $i$ ($0 \le i \le d$).
Let ${\cal A}$ have minimal idempotents $E_i$ ($0 \le i \le d$).
The column spaces of the $E_i$ are common eigenspaces of all
matrices in ${\cal A}$. Let $P_{ij}$ be the corresponding eigenvalue
of $A_j$, so that $A_jE_i = P_{ij}E_i$ ($0 \le i,j \le d$).
Now $A$ has eigenvalues $\theta_i = P_{i1}$ with multiplicities
$m_i = \rk E_i = \tr E_i$ ($0 \le i \le d$).
Index the eigenvalues such that $\theta_0 > \theta_1 > \cdots > \theta_d$.

\medskip
Standard facts about Sturm sequences give information on
the sign pattern of the matrix $P$.

\begin{Proposition}
Let $\Gamma$ be distance-regular, and $P$ its eigenvalue matrix.
Then row $i$ and column $i$ $(0\le i\le d)$ of $P$ both have $i$ sign changes.
In particular, row $d$ and column $d$ consist of nonzero numbers
that alternate in sign. \qed
\end{Proposition}

If $M \in {\cal A}$ and $0 \le i \le d$, then
$M \prod_{j \ne i} (A - \theta_j I) = c(M,i) E_i$
for some constant $c(M,i)$.
We apply this observation to $M = A_d$.

\begin{Proposition}
Let $\Gamma$ have intersection array
$\{b_0,b_1,\ldots,b_{d-1};\,c_1,c_2,\ldots,c_d\}$. Then for each $i$
$(0 \le i \le d)$ we have
\begin{equation}\label{eq1}
m_i A_d \prod_{j \ne i} (A - \theta_j I) = nb_0b_1\cdots b_{d-1} E_i .
\end{equation}
\end{Proposition}
\Proof
Both sides differ by a constant factor. Take traces on both sides.
Since $\tr A_d A^h = 0$ for $h < d$ it follows that
$\tr A_d \prod_{j \ne i} (A - \theta_j I) = \tr A_d A^d =
c_1c_2\cdots c_dnk_d = nb_0b_1\cdots b_{d-1}$.
Now the result follows from $\tr E_i = m_i$.
\qed

This can be said in an equivalent numerical way.

\begin{Corollary}
We have
~$m_i P_{id} \prod_{j \ne i} (\theta_i - \theta_j) = nb_0b_1\cdots b_{d-1}$~
for each $i$ $(0\le i\le d)$.
\end{Corollary}
\Proof Multiply \eqref{eq1} by $E_i$. \qed

We find a criterion for $A_d$ to have
two equal eigenvalues $P_{gd}$ and $P_{hd}$.

\begin{Proposition}\label{eq2}
For $g \ne h$,
$P_{gd} = P_{hd}$ if and only if
$\sum_i m_i \prod_{j \ne g,h} (\theta_i - \theta_j) = 0$.
\end{Proposition}
\Proof
$P_{gd} = P_{hd}$ if and only if
$m_g \prod_{j \ne g} (\theta_g - \theta_j) =
m_h \prod_{j \ne h} (\theta_h - \theta_j)$. \qed

For example, the Biggs-Smith graph has diameter $d=7$
and spectrum $3^1$, $\theta_1^9$, $2^{18}$, $\theta_3^{16}$,
$0^{17}$, $\theta_5^{16}$, $\theta_6^9$, $\theta_7^{16}$,
where $\theta_i$, $i=1,6$, satisfy $f(\theta) = \theta^2-\theta-4 = 0$
and $\theta_i$, $i=3,5,7$, satisfy $g(\theta) = \theta^3+3\theta^2-3 = 0$.
Now $P_{27} = P_{47}$ since
$$\sum_i m_i \prod_{j \ne 2,4} (\theta_i - \theta_j) =
\sum_{i=2,4} m_i (\theta_i-3)f(\theta_i)g(\theta_i) = 0.$$

One can generalize Proposition \ref{eq2}, and see:
\begin{Proposition} \label{multipleeq}
Let $H \subseteq \{0,\ldots,d\}$. Then all $P_{hd}$ for $h \in H$
take the same value if and only if
$\sum_i m_i \theta_i^e \prod_{j \notin H} (\theta_i - \theta_j) = 0$
for $0 \le e \le |H|-2$.
\end{Proposition}
\Proof
Induction on $|H|$. We just did the case $|H|=2$.
Let $|H| > 2$ and let $h,h' \in H$. We do the `only if' part.
By induction $\sum_i m_i \theta_i^e \prod_{j \notin H\setminus\{x\}}
(\theta_i - \theta_j) = 0$ holds for $0 \le e \le |H|-3$
and $x = h,h'$. Subtract these two formulas and divide by
$\theta_h - \theta_{h'}$ to get
$\sum_i m_i \theta_i^e \prod_{j \notin H}
(\theta_i - \theta_j) = 0$ for $0 \le e \le |H|-3$.
Then add the first formula for $x=h$ and $\theta_h$ times the last
formula, to get the same conclusion for $1 \le e \le |H|-2$.
The converse is clear.
\qed

Since the $P_{id}$ alternate in sign, the largest sets $H$ that can
occur here are $\{0,2,\ldots,d\}$ and $\{1,3,\ldots,d-1\}$ for $d=2e$
and $\{0,2,\ldots,d-1\}$ and $\{1,3,\ldots,d\}$ for $d=2e+1$.
We investigate such sets below 
(see `the half-antipodal case').

\medskip
For small $d$ one can use identities like
$\sum_i m_i = n$, $\sum_i m_i \theta_i = 0$,
$\sum_i m_i \theta_i^2 = nk$,
$\sum_i m_i \theta_i^3 = nk\lambda$ (where $k=b_0$,
and $\lambda = k-1-b_1$) to simplify the condition
of Proposition \ref{eq2}.
Let us do some examples. Note that $\theta_0 = k$.

\subsection*{The case \boldmath$d=3$}
For $d=3$ we find that $P_{13} = P_{33}$ if and only if
$\sum_i m_i (\theta_i-\theta_0)(\theta_i-\theta_2) = 0$,
i.e., if and only if $nk+n\theta_0\theta_2 = 0$,
i.e., if and only if $\theta_2 = -1$
(cf.~\cite{BCN}, 4.2.17).

\subsection*{The case \boldmath$d=4$}
For $d=4$ we find that $P_{14} = P_{34}$ if and only if
$\sum_i m_i (\theta_i-\theta_0)(\theta_i-\theta_2)(\theta_i-\theta_4) = 0$,
i.e., if and only if
$nk\lambda - nk(\theta_0+\theta_2+\theta_4)-n\theta_0\theta_2\theta_4 = 0$.
This happens if and only if
$(\theta_2+1)(\theta_4+1) = -b_1$.
Of course $P_{24} = P_{44}$ will follow from
$(\theta_1+1)(\theta_3+1) = -b_1$.

\medskip\noindent
A generalized octagon ${\rm GO}(s,t)$ has eigenvalues
$$
(\theta_i)_i = (s(t+1),~ s-1+\sqrt{2st},~ s-1,~
s-1-\sqrt{2st},~ -t-1)
$$
and $b_1 = st$.
Since $(\theta_2+1)(\theta_4+1) = -b_1$ it follows that
$P_{14} = P_{34}$, so that $\Gamma_4$ does not have
more than 4 distinct eigenvalues.

\medskip\noindent
A dual polar graph ${}^2D_5(q)$ has eigenvalues
$$(\theta_i)_i = (q^5{+}q^4{+}q^3{+}q^2,~
q^4{+}q^3{+}q^2{-}1,~ q^3{+}q^2{-}q{-}1,~ {-}q{-}1,~ {-}q^3{-}q^2{-}q{-}1)$$
and $b_1 = q^3 (q^2+q+1)$.
Since $(\theta_1+1)(\theta_3+1) = -b_1$ it follows that
$P_{24} = P_{44}$.

\medskip\noindent
Some further examples:

\medskip\noindent
\begin{tabular}{@{~}l@{~}c@{~~}l@{~~}l@{~~}l@{}}
name & $n$ & intersection array & spectrum & equality \\
\hline
Coxeter graph & 28 & $\{3,2,2,1;\,1,1,1,2\}$ &
$3^1$ $2^8$ $a^6$ $(-1)^7$ $b^6$ &
$P_{14} = P_{34}$ \\
&&&$a,b=-1\pm\sqrt{2}$ \\
Odd graph $O_5$ & 126 & $\{5,4,4,3;\,1,1,2,2\}$ &
$5^1$ $3^{27}$ $1^{42}$ $(-2)^{48}$ $(-4)^8$ &
$P_{24} = P_{44}$ \\
$M_{22}$ graph & 330 & $\{7,6,4,4;\,1,1,1,6\}$ &
$7^1$ $4^{55}$ $1^{154}$ $(-3)^{99}$ $(-4)^{21}$ &
$P_{14} = P_{34}$ \\
Unital graph & 280 & $\{9,8,6,3;\,1,1,3,8\}$ &
$9^1$ $4^{64}$ $1^{105}$ $(-3)^{90}$ $(-5)^{20}$ &
$P_{14} = P_{34}$ \\
\end{tabular}

\subsection*{The case \boldmath$d=4$ with strongly regular $\Gamma_4$}

One may wonder whether it is possible that $\Gamma_4$ is
strongly regular. This would require $p^1_{44} = p^2_{44} = p^3_{44}$.
Or, equivalently, that $\Gamma_4$ has only two eigenvalues
with eigenvector other than the all-1 vector. Since the
values $P_{i4}$ alternate in sign, this would mean
$P_{14} = P_{34}$ and $P_{24} = P_{44}$.

\begin{Proposition}
Let $\Gamma$ be a distance-regular graph of diameter $4$.
The following assertions are equivalent.

(i) $\Gamma_4$ is strongly regular.

(ii) $b_3 = a_4+1$ and $b_1=b_3c_3$.

(iii) $(\theta_1+1)(\theta_3+1) = (\theta_2+1)(\theta_4+1) = -b_1$.
\end{Proposition}
%

\medskip\noindent
\Proof
(i)-(ii)
A boring computation (using \cite{BCN}, 4.1.7) shows that
$p^1_{44} = p^2_{44}$ is equivalent to $b_3 = a_4+1$, and that
if this holds $p^1_{44} = p^3_{44}$ is equivalent to $b_1=b_3c_3$.

(i)-(iii)
$\Gamma_4$ will be strongly regular if and only if
$P_{14} = P_{34}$ and $P_{24} = P_{44}$.
We saw that this is equivalent to
$(\theta_2+1)(\theta_4+1) = -b_1$ and $(\theta_1+1)(\theta_3+1) = -b_1$.
\qed

\noindent
The fact that (i) implies the first equality in (iii) was proved in
\cite{F01} as a consequence of another characterization of (i)
in terms of the spectrum only.
More generally, a quasi-spectral characterization of those connected
regular graphs (with $d+1$ distinct eigenvalues) which are distance-regular,
and with the distance-$d$ graph being strongly regular, is given in
\cite[Th. 2.2]{F00}.

\medskip
No nonantipodal examples are known, but the infeasible array
$\{12,8,6,4;\,1,\maysplit 1,\maysplit 2,9\}$ with spectrum
$12^1$ $7^{56}$ $3^{140}$ $(-2)^{160}$ $(-3)^{168}$
(cf.~\cite{BCN},\,p.\,410) would~have been an example
(and there are several open candidate arrays, such as
$\{21,20,\maysplit 14,\maysplit 10;\,1,1,2,12\}$,
$\{24,20,20,10;\,1,1,2,15\}$, and
$\{66,65,63,13;\,1,1,5,54\}$).
%
%

\medskip
If $\Gamma$ is antipodal, then $\Gamma_4$ is a union of cliques
(and hence strongly regular). This holds precisely when
$\theta_1 + \theta_3 = \lambda$ and $\theta_1\theta_3 = -k$
and $(\theta_2+1)(\theta_4+1) = -b_1$.
(Indeed, $\theta_1,\theta_3$ are the two roots of
$\theta^2-\lambda\theta-k = 0$ by \cite{BCN}, 4.2.5.).  \\
There are many examples, e.g.

\medskip\noindent
\begin{tabular}{@{~}l@{~}c@{~~}l@{~~}l@{}}
name & $n$ & intersection array & spectrum \\
\hline
Wells graph & 32 & $\{5,4,1,1;\,1,1,4,5\}$ &
$5^1$ $\sqrt{5}{}^8$ $1^{10}$ $(-\sqrt{5})^8$ $(-3)^5$ \\
3.Sym(6).2 graph & 45 & $\{6,4,2,1;\,1,1,4,6\}$ &
$6^1$ $3^{12}$ $1^9$ $(-2)^{18}$ $(-3)^5$ \\
Locally Petersen & 63 & $\{10,6,4,1;\,1,2,6,10\}$ &
$10^1$ $5^{12}$ $1^{14}$ $(-2)^{30}$ $(-4)^6$
\end{tabular}

\medskip
If $\Gamma$ is bipartite, then $\Gamma_4$ is disconnected,
so if it is strongly regular, it is a union of cliques
and $\Gamma$ is antipodal. In this case its spectrum is
$$\{k^1,~\sqrt{k}{\,}^{n/2-k},~
0^{2k-2},~ (-\sqrt{k})^{n/2-k},~ (-k)^1 \}.$$
Such graphs are precisely the incidence graphs of
symmetric $(m,\mu)$-nets, where $m = k/\mu$ (\cite{BCN}, p. 425).

\subsection*{The case \boldmath$d=5$}
As before, and also using
$\sum_i m_i \theta_i^4 = nk(k+\lambda^2+b_1\mu)$
(where $\mu = c_2$) we find for $\{f,g,h,i,j\} = \{1,2,3,4,5\}$
that $P_{fd} = P_{gd}$ if and only if
$$(\theta_h+1)(\theta_i+1)(\theta_j+1) + b_1(\theta_h+\theta_i+\theta_j) =
b_1(\lambda-\mu-1).$$

In case $\theta_i = -1$,
this says that $\theta_h+\theta_j = \lambda-\mu$.

For example, the Odd graph $O_6$ has $\lambda=0$, $\mu=1$,
and eigenvalues 6, 4, 2, $-1$, $-3$, $-5$.
It follows that $P_{15} = P_{55}$ and $P_{25} = P_{45}$.

Similarly, the folded 11-cube has $\lambda=0$, $\mu=2$,
and eigenvalues 11, 7, 3, $-1$, $-5$, $-9$.
It follows that $P_{15} = P_{55}$ and $P_{25} = P_{45}$.

An example without eigenvalue $-1$ is provided by the
folded Johnson graph $\bar{J}(20,10)$. It has $\lambda=18$,
$\mu=4$ and eigenvalues 100, 62, 32, 10, $-4$, $-10$.
We see that $P_{35} = P_{55}$.

\medskip\noindent
Combining two of the above conditions, we see that
$P_{15} = P_{35} = P_{55}$ if and only if
$(\theta_2+1)(\theta_4+1) = -b_1$ and $\theta_2 + \theta_4 = \lambda - \mu$
(and hence $\theta_2\theta_4 = \mu-k$).
%
\\ Now $b_3+b_4+c_4+c_5 = 2k+\mu-\lambda$ and
$b_3b_4+b_3c_5+c_4c_5 = kb_1+k\mu+\mu$.

\subsection*{Generalized 12-gons}
A generalized 12-gon of order $(q,1)$ (the line graph
of the bipartite point-line incidence graph
of a generalized hexagon of order $(q,q)$) has diameter 6,
and its $P$ matrix is given by

{\small
$$
P \!=\! \left(\!\begin{array}{ccccccc}
1 & 2q & 2q^2 & 2q^3 & 2q^4 & 2q^5 & q^6 \\
1 & q{-}1{+}a & q{+}(q{-}1)a & 2q(q{-}1) & {-}q^2{+}q(q{-}1)a & q^2(q{-}1){-}q^2a & {-}q^3 \\
1 & q{-}1{+}b & {-}q{+}(q{-}1)b & {-}2qb & {-}q^2{-}q(q{-}1)b & {-}q^2(q{-}1){+}q^2b & q^3 \\
1 & q{-}1 & {-}2q & {-}q(q{-}1) & 2q^2 & q^2(q{-}1) & {-}q^3 \\
1 & q{-}1{-}b & {-}q{-}(q{-}1)b & 2qb & {-}q^2{+}q(q{-}1)b & {-}q^2(q{-}1){-}q^2b & q^3 \\
1 & q{-}1{-}a & q{-}(q{-}1)a & 2q(q{-}1) & {-}q^2{-}q(q{-}1)a & q^2(q{-}1){+}q^2a & {-}q^3 \\
1 & {-}2 & 2 & {-}2 & 2 & {-}2 & 1
\end{array}\!\right)
$$

}\noindent
where $a = \sqrt{3q}$ and $b = \sqrt{q}$.
We see that $P_{16} = P_{36} = P_{56}$ and $P_{26} = P_{46}$.

\medskip
Its dual is a generalized 12-gon of order $(1,q)$, and is bipartite.
The $P$ matrix is given by

{\small
$$
P \!=\! \left(\!\begin{array}{ccccccc}
1 & q+1 & q(q+1) & q^2(q+1) & q^3(q+1) & q^4(q+1) & q^5 \\
1 & a & 2q-1 & (q-1)a & q(q-2) & -qa & -q^2 \\
1 & b & -1 & -(q+1)b & -q^2 & qb & q^2 \\
1 & 0 & -q-1 & 0 & q(q+1) & 0 & -q^2 \\
1 & -b & -1 & (q+1)b & -q^2 & -qb & q^2 \\
1 & -a & 2q-1 & -(q-1)a & q(q-2) & qa & -q^2 \\
1 & -q-1 & q(q+1) & -q^2(q+1) & q^3(q+1) & -q^4(q+1) & q^5
\end{array}\!\right)
$$

}\noindent
where $a = \sqrt{3q}$ and $b = \sqrt{q}$.
We see that $P_{16} = P_{36} = P_{56}$ and $P_{26} = P_{46}$.\\
As expected (cf.~\cite{B}), the squares of all $P_{id}$ are powers of $q$.

\subsection*{Dual polar graphs}
According to \cite{B}, dual polar graphs of diameter $d$
satisfy $$P_{id} = (-1)^i q^{d(d-1)/2+de-i(d+e-i)}$$
where $e$ has the same meaning as in \cite{BCN}, 9.4.1.
It follows that $P_{hd} = P_{id}$ when $d+e$ is even
and $h+i = d+e$.

\medskip\noindent
For the dual polar graphs $B_d(q)$ and $C_d(q)$ we have $e=1$,
and the condition becomes $h+i = d+1$ where $d$ is odd.
Below we will see this in a different way.

For the dual polar graph $D_d(q)$ we have $e = 0$,
and the condition becomes $h+i = d$ where $d$ is even.
Not surprising, since this graph is bipartite.

For the dual polar graph ${}^2D_{d+1}(q)$ we have $e=2$,
and the condition becomes $h+i = d+2$ where $d$ is even.
(We saw the case $d=4$ above.)

Finally, $h+i = d+e$ is impossible when $e$ is not integral.

\subsection*{Distance-regular distance 1-or-2 graph}
The distance 1-or-2 graph $\Delta = \Gamma_1 \cup \Gamma_2$ of $\Gamma$
(with adjacency matrix $A_1+A_2$) is distance-regular
if and only if $b_{i-1}+b_i+c_i+c_{i+1} = 2k+\mu-\lambda$
for $1 \le i \le d-1$, cf.~\cite{BCN}, 4.2.18.

\begin{Proposition}
Suppose that $\Gamma_1 \cup \Gamma_2$ is distance-regular.
Then for $1 \le i \le d$ we have
$P_{d+1-i,d} = P_{id}$ if $d$ is odd,
and $(\theta_{d+1-i}+1)P_{i,d} =
(\theta_i+1)P_{d+1-i,d}$ if $d$ is even.
If $i \ne (d+1)/2$ then $\theta_{d+1-i} = \lambda-\mu-\theta_i$.
If $d$ is odd, then $\theta_{(d+1)/2} = -1$.
\end{Proposition}
\Proof
For each eigenvalue $\theta$ of $\Gamma$, there is
an eigenvalue $(\theta^2 + (\mu-\lambda)\theta - k)/\mu$ of $\Delta$.
If $d$ is odd, then $\Delta$ has diameter $(d+1)/2$,
and $\Gamma$ has an eigenvalue $-1$,
and for each eigenvalue $\theta \ne k, -1$ of $\Gamma$
also $\lambda-\mu-\theta$ is an eigenvalue.
Now $\Delta_{(d+1)/2} = \Gamma_d$, and $P_{i'd} = P_{id}$
if $i,i'$ belong to the same eigenspace of $\Delta$.
Since the numbers $P_{id}$ alternate, the eigenvalue $-1$
must be the middle one (not considering $\theta_0$), and we see that
$P_{d+1-i,d} = P_{id}$ for $1 \le i \le d$.

If $d$ is even, then $\Delta$ has diameter $d/2$, and for each
eigenvalue $\theta \ne k$ of $\Gamma$ also $\lambda-\mu-\theta$
is an eigenvalue. Now $\Delta_{d/2} = \Gamma_{d-1} \cup \Gamma_d$.
Since $AA_d = b_{d-1}A_{d-1} + a_dA_d$ and our parameter conditions
imply $b_{d-1}+c_d = b_1+\mu$
(so that $b_{d-1}(P_{i,d-1} + P_{i,d}) = (\theta_i - a_d + b_{d-1})P_{i,d} =
(\theta_i+\mu-\lambda-1)P_{i,d}$), the equalities
$P_{i,d-1} + P_{i,d} = P_{d+1-i,d-1} + P_{d+1-i,d}$ $(1 \le i \le d)$
and $\theta_i + \theta_{d+1-i} = \lambda-\mu$
imply $(\theta_{d+1-i}+1)P_{i,d} = (\theta_i+1)P_{d+1-i,d}$.
\qed

\medskip
The Odd graph $O_{d+1}$ on $\binom{2d+1}{d}$ vertices
has diameter $d$ and eigenvalues $\theta_i = d+1-2i$ for $i < (d+1)/2$,
and $\theta_i = d-2i$ for $i \ge (d+1)/2$.
Since its distance 1-or-2 graph is distance-regular,
we have $P_{d+1-i,d} = P_{id}$ for odd $d$ and $1 \le i \le d$.
%

The folded $(2d+1)$-cube on $2^{2d}$ vertices
has diameter $d$ and eigenvalues
$\theta_i = 2d+1-4i$ with multiplicities
$m_i = \binom{2d+1}{2i}$ ($0 \le i \le d$).
Since its distance 1-or-2 graph is distance-regular,
it satisfies $P_{d+1-i,d} = (-1)^{d+1} P_{id}$ for $1 \le i \le d$.
(Note that $\theta_{d+1-i}+1 = -(\theta_i+1)$ since $\mu-\lambda = 2$.)

The dual polar graphs $B_d(q)$ and $C_d(q)$ have diameter $d$
and eigenvalues $\theta_i = (q^{d-i+1}-q^i)/(q-1) - 1$.
Since their distance 1-or-2 graphs are distance-regular,
they satisfy $P_{d+1-i,d} = (-1)^{d+1} P_{id}$ for $1 \le i \le d$.
%

\medskip
The fact that $-1$ must be the middle eigenvalue for odd $d$,
implies
that $\theta_{(d-1)/2}> \lambda-\mu+1$,
so that there is no eigenvalue $\xi$
with $-1 < \xi < \lambda-\mu+1$.

\subsection*{The bipartite case}
If $\Gamma$ is bipartite, then $\theta_{d-i} = -\theta_i$,
and $P_{d-i,j} = (-1)^j P_{i,j}$ ($0 \le i,j \le d$).
In particular, if $d$ is even, then $P_{d-i,d} = P_{id}$
and $\Gamma_d$ is disconnected.

\subsection*{The antipodal case}
The graph $\Gamma$ is antipodal when having distance $d$ is an
equivalence relation, i.e., when $\Gamma_d$ is a union of cliques.
The graph is called an antipodal $r$-cover, when these cliques
are $r$-cliques. Now $r = k_d+1$, and
$P_{id}$ alternates between $k_d$ and $-1$.

For example, the ternary Golay code graph (of diameter 5)
with intersection array
$\{22,20,18,2,1;\,1,2,9,20,22\}$ has spectrum
$22^1$ $7^{132}$ $4^{132}$ $(-2)^{330}$ $(-5)^{110}$ $(-11)^{24}$
and satisfies $P_{05} = P_{25} = P_{45} = 2$,
$P_{15} = P_{35} = P_{55} = -1$.

For an antipodal distance-regular graph $\Gamma$, the folded graph has
eigenvalues $\theta_0,\theta_2,\ldots,\theta_{2e}$ where $e = [d/2]$.
In Theorem \ref{oddhalfatp} below we show for odd $d$ that this
already follows from $P_{1d} = P_{3d} = \cdots = P_{dd}$.

\begin{Proposition}
\label{0ddd}
If $P_{0d} = P_{id}$ then $i$ is even.
Let $i > 0$ be even. Then $P_{0d} = P_{id}$ if and only
$\Gamma$ is antipodal, or $i=d$ and $\Gamma$ is bipartite.
\end{Proposition}
\Proof
Since the $P_{id}$ alternate in sign, $P_{0d} = P_{id}$ implies
that $i$ is even.
If $\Gamma$ is bipartite, then $P_{dd} = (-1)^d P_{0d}$.
If $\Gamma$ is antipodal, then $P_{id} = P_{0d}$ for all even $i$.
That shows the `if' part. Conversely, if $P_{0d} = P_{id}$,
then the valency of $\Gamma_d$ is an eigenvalue
of multiplicity larger than 1, so that $\Gamma_d$ is disconnected,
and hence $\Gamma$ is imprimitive and therefore antipodal or bipartite.
If $\Gamma$ is bipartite but not antipodal, then its halved graphs
are primitive and $|P_{id}| < P_{0d}$ for $0 < i < d$
(cf.~\cite{BCN}, pp.~140--141).
\qed

\subsection*{The half-antipodal case} 
Given an array $\{b_0,\ldots,b_{d-1};\,c_1,\ldots,c_d\}$
of positive real numbers, define the polynomials $p_i(x)$ for
$-1 \le i \le d+1$ by $p_{-1}(x)=0$, $p_0(x)=1$,
$(x-a_i)p_i(x) = b_{i-1}p_{i-1}(x)+c_{i+1}p_{i+1}(x)$
($0 \le i \le d$), where $a_i = b_0-b_i-c_i$ and $c_{d+1}$
is some arbitrary positive number.
The eigenvalues of the array are by definition the zeros of
$p_{d+1}(x)$, and do not depend on the choice of $c_{d+1}$.
Each $p_i(x)$ has degree $i$, and, by the theory of Sturm sequences,
each $p_i(x)$ has $i$ distinct real zeros, where the zeros of
$p_{i+1}(x)$ interlace those of $p_i(x)$.
If $\{b_0,\ldots,b_{d-1};\,c_1,\ldots,c_d\}$ is the intersection array
of a distance-regular graph $\Gamma$, then the eigenvalues of the array
are the eigenvalues of (the adjacency matrix of) $\Gamma$.

Let $L$ be the tridiagonal matrix
$$
L=\left(\begin{array}{cccccc}
a_0 & b_0 \\
c_1 & a_1 & b_1 \\
  & c_2 & . & . \\
  &     & . & . & . \\
  &     &   & . & . & b_{d-1} \\
  &     &   &   & c_d & a_d
\end{array}\right).
$$
The eigenvalues of the array $\{b_0,\ldots,b_{d-1};\,c_1,\ldots,c_d\}$
are the eigenvalues of the matrix $L$.

\begin{Theorem} \label{oddhalfatp}
Let $\Gamma$ be a distance-regular graph with odd diameter $d = 2e+1$
and intersection array $\{b_0,\ldots,b_{d-1};\,c_1,\ldots,c_d\}$.
Then $P_{1d} = P_{3d} = \cdots = P_{dd}$ if and only if
the $\theta_j$ with $j=0,2,4,\ldots,2e$ are the eigenvalues of the array
$\{b_0,\ldots,b_{e-1};\,c_1,\ldots,c_e\}$.
\end{Theorem}
\Proof
Let $H = \{1,3,\ldots,d\}$, so that $|H| = e+1$.
By Proposition \ref{multipleeq}, $P_{1d} = P_{3d} = \cdots = P_{dd}$
if and only if
$\sum_i m_i \theta_i^s \prod_{j \notin H} (\theta_i - \theta_j) = 0$
for $0 \le s \le e-1$.
Let~$E = \{0,2,\ldots,2e\}$, so that $|E| = e+1$.
Then this condition is equivalent to
$$
\tr A^s \prod_{j \in E} (A-\theta_j I) = 0\qquad (0 \le s \le e-1).
$$
This says that the expansion of $\prod_{j \in E} (A-\theta_j I)$
in terms of the $A_i$ does not contain $A_s$ for $0 \le s \le e-1$,
hence is equivalent to
$\prod_{j \in E} (A-\theta_j I) = aA_e + bA_{e+1}$ for certain constants $a,b$.
Since $0 \in E$, we find that $ak_e+bk_{e+1} = 0$, and the condition is
equivalent to $(A_e/k_e - A_{e+1}/k_{e+1})E_j = 0$ for all $j \in E$.

An eigenvalue $\theta$ of $\Gamma$ defines a right eigenvector $u$
(known as the `standard sequence') by $Lu = \theta u$.
It follows that $\theta$ will be an eigenvalue of the array
$\{b_0,\ldots,b_{e-1};\,c_1,\ldots,c_e\}$ precisely when $u_e = u_{e+1}$.
Up to scaling, the $u_i$ belonging to $\theta_j$ are the $Q_{ij}$
(that is, the columns of $Q$ are eigenvectors of $L$).
So, $\theta_j$ is an eigenvalue of $\{b_0,\ldots,b_{e-1};\,c_1,\ldots,c_e\}$
for all $j \in E$ precisely when $Q_{ej} = Q_{e+1,j}$ for all $j \in E$.
Since $k_iQ_{ij} = m_jP_{ji}$, this holds if and only if
$P_{je}/k_e = P_{j,e+1}/k_{e+1}$, i.e., if and only if
$(A_e/k_e - A_{e+1}/k_{e+1})E_j = 0$ for all $j \in E$.
\qed

For example, if $d=3$ one has $P_{13} = P_{33}$ if and only if
$\theta_0,\theta_2$ are the eigenvalues $k$, $-1$ of
the array $\{k;\,1\}$. And if $d=5$ one has
$P_{15} = P_{35} = P_{55}$ if and only if
$\theta_0, \theta_2, \theta_4$ are the eigenvalues of the array
$\{k,b_1;\,1,c_2\}$.

\medskip\noindent
The case of even $d$ is slightly more complicated.

\begin{Theorem}
Let $\Gamma$ be a distance-regular graph with even diameter $d = 2e$
and intersection array $\{b_0,\ldots,b_{d-1};\,c_1,\ldots,c_d\}$.
Then $P_{1d} = P_{3d} = \cdots = P_{d-1,d}$ if and only if
the $\theta_j$ with $j=0,2,4,\ldots,2e$ are the eigenvalues of the array
$\{b_0,\ldots,b_{e-1};\,c_1,\ldots,c_{e-1},c_e+zb_e\}$ for some
real number $z$ with $0 < z \le 1$, uniquely determined by
$\sum_{i=0}^e \theta_{2i} = \sum_{i=0}^e a_i + (1-z)b_e$.
If \,$\Gamma$ is antipodal or bipartite, then $z=1$.
\end{Theorem}
\Proof
Let $E = \{0,2,\ldots,d\}$.
As before we see that $P_{1d} = P_{3d} = \cdots = P_{d-1,d}$
is equivalent to the condition that
$\prod_{j \in E} (A-\theta_j I) = aA_{e-1} + bA_e + cA_{e+1}$
for certain constants $a,b,c$.
Comparing coefficients of $A^{e+1}$ we see that $c > 0$.
With $j=0$ we see that $ak_{e-1}+bk_e+ck_{e+1} = 0$.
%

Take
$$
z = -\frac{ak_{e-1}}{ck_{e+1}} = 1 + \frac{bk_e}{ck_{e+1}}.$$
Then $aP_{j,e-1}+bP_{je}+cP_{j,e+1} = 0$ for $j \in E$ gives
$$
z\left(\frac{P_{j,e-1}}{k_{e-1}} - \frac{P_{je}}{k_e}\right)
 = \frac{P_{j,e+1}}{k_{e+1}} - \frac{P_{je}}{k_e}.
$$

On the other hand, if $\theta = \theta_j$ for some $j \in E$, then
$\theta$ is an eigenvalue of the array
$\{b_0,\ldots,b_{e-1};\,c_1,\ldots,c_{e-1},\maysplit c_e{+}zb_e\}$
precisely when $c_eu_{e-1}+(k{-}b_e{-}c_e)u_e+b_eu_{e+1} =
(c_e{+}zb_e)u_{e-1}+(k{-}c_e{-}zb_e)u_e$,
i.e., when $z(u_{e-1}-u_e) = u_{e+1}-u_e$.
Since (up to a constant factor) $u_i = P_{ji}/k_i$, this is equivalent
to the condition above.

Noting that $d \in E$, we can apply the above to $\theta = \theta_d$.
Since the bottom row of $P$ has $d$ sign changes, it follows that
the sequence $u_i$ has $d$ sign changes.
In particular, the $u_i$ are nonzero.
Now $u_{e-1}-u_e$ and $u_{e+1}-u_e$ have the same sign,
and it follows that $z > 0$.

If $\Gamma$ is an antipodal $r$-cover of diameter $d = 2e$, then
$c_e+zb_e = rc_e$ and $z = 1$ (and $P_{je} = 0$ for all odd $j$).

If $\Gamma$ is bipartite, then $\theta_i + \theta_{d-i} = 0$
for all $i$, so $\sum_{j \in E} \theta_j = 0$, so our tridiagonal
matrix (the analog of $L$) has trace
$0 = a_1+\cdots+a_e+(1-z)b_e = (1-z)b_e$, so that $z = 1$.

It remains to show that $z \le 1$. We use $z(u_{e-1}-u_e) = u_{e+1}-u_e$
to conclude that $\theta$ is an eigenvalue of both
$$
\left(\begin{array}{c@{~~}c@{~~}c@{~~}c@{~~}c}
0 & k \\
c_1 & a_1 & b_1 \\
  & . & . & . \\
  &   & c_{e-1} & a_{e-1} & b_{e-1} \\
  &   &   & p & k-p
\end{array}\right)
~{\rm and}~
\left(\begin{array}{c@{~~}c@{~~}c@{~~}c@{~~}c}
k-q & q \\
c_{e+1} & a_{e+1} & b_{e+1} \\
  &  . & . & . \\
  &     & c_{d-1} & a_{d-1} & b_{d-1} \\
  &     &   & c_d & a_d
\end{array}\right),
$$
where $p = c_e+zb_e$ and $q = b_e+z^{-1}c_e$.
Since this holds for each $\theta=\theta_j$ for $j \in E$,
this accounts for all eigenvalues of these two matrices, and
$\sum_{j \in E} \theta_j = a_1+\cdots+a_e+(1-z)b_e =
a_e+\cdots+a_d+(1-z^{-1})c_e$.
Since $p^{e+1}_{ed} \ge 0$ it follows that
$a_1+\cdots+a_{e-1} \le a_{e+1}+\cdots+a_d$ (cf.~\cite{BCN}, 4.1.7),
and therefore $(1-z)b_e \ge (1-z^{-1})c_e$.
It follows that $z \le 1$.
\qed

%
%

\medskip

% ...

\medskip
Nonantipodal, nonbipartite examples:

\medskip\noindent
\begin{tabular}{@{\,}llll@{\,}}
name & array & half array & $z$ \\
\hline
Coxeter & $\{3,2,2,1;\,1,1,1,2\}$ & $\{3,2;\,1,2\}$ & $1/2$ \\
$M_{22}$ & $\{7,6,4,4;\,1,1,1,6\}$ & $\{7,6;\,1,3\}$ & $1/2$ \\
$P\Gamma L(3,4).2$ & $\{9,8,6,3;\,1,1,3,8\}$ & $\{9,8;\,1,4\}$ & $1/2$ \\
gen.\,8-gon & $\{s(t+1),st,st,st;\,1,1,1,t+1\}$
 & $\{s(t+1),st;\,1,t+1\}$ & $1/s$ \\
gen.\,12-gon & $\{2q,q,q,q,q,q;\,1,1,1,1,1,2\}$
 & $\{2q,q,q;\,1,1,2\}$ & $1/q$
\end{tabular}

\medskip
Concerning the value of $z$, note that both $zb_e$ and $z^{-1}c_e$
are algebraic integers.

\bigskip
If $d=4$, the case $z = 1$ can be classified.

\begin{Proposition}
Let $\Gamma$ be a distance-regular graph with even diameter $d = 2e$
such that $\theta_j$ with $j=0,2,4,\ldots,2e$ are the eigenvalues
of the array $\{b_0,\ldots,b_{e-1};\maysplit\,c_1,\ldots,c_{e-1},b_e+c_e\}$.
Then $\Gamma$ satisfies $p^d_{e,e+1} = 0$.
If moreover $d \le 4$, then $\Gamma$ is antipodal or bipartite.
\end{Proposition}
\Proof
We have equality in the inequality $a_1+\cdots+a_{e-1} \le
a_{e+1}+\cdots+a_d$, so that $p^{e+1}_{ed} = 0$ by \cite{BCN}, 4.1.7.
%
The case $d=2$ is trivial. Suppose \mbox{$d=4$}. Then $p^3_{24} = 0$.
Let $d(x,z) = 4$, and consider neighbors $y,w$ of $z$,
where $d(x,y) = 3$ and $d(x,w) = 4$. Then $d(y,w) \ne 2$ since
there are no 2-3-4 triangles, so $d(y,w) = 1$.
If $a_4 \ne 0$ then there exist such vertices $w$, and we
find that the neighborhood $\Gamma(z)$ of $z$ in $\Gamma$
is not coconnected (its complement is not connected),
contradicting \cite{BCN}, 1.1.7. Hence $a_4 = 0$, and $a_3 = a_1$.
If $b_3 > 1$, then let $d(x,y)=3$, $y \adj z,z'$ with $d(x,z)=d(x,z')=4$.
Let $z''$ be a neighbor of $z'$ with $d(z,z'') = 2$.
Then $d(x,z'')=3$ since $a_4=0$, and we see a 2-3-4 triangle, contradiction.
So if $b_3 > 1$ then $a_2 = 0$, and $a_1 = 0$ since $a_1 \le 2a_2$,
and the graph is bipartite. If $b_3=1$, then the graph is antipodal.
\qed

\subsection*{Variations}
One can vary the above theme. First, the following general result builds on ideas previously used.

\begin{Theorem}
Let $\Gamma$ be a distance-regular graph with diameter $d$, and
let $H \subseteq \{0,\ldots,d\}$. Then all $P_{id}$ for $i \in H$
take the same value  if and only if
the eigenvalues $\theta_j$ with $j\notin H$ are the zeros of the polynomial
$\sum_{i=|H|-1}^{r} \frac{\alpha_i}{k_i}p_i(x)$, where  $r=d+1-|H|$,
$\alpha_{r}=1$, and
\begin{equation}
\label{eqThm}
\alpha_i =\frac{\tr (A_i \prod_{j\notin H}(A-\theta_j I))}{\tr (A_r A^{r})}\qquad
(|H|-1\le i\le d-|H|).
\end{equation}
\end{Theorem}
\Proof
By Proposition \ref{multipleeq}, the hypothesis holds if and only if we have the equalities
$\sum_i m_i \theta_i^s \prod_{j \notin H} (\theta_i - \theta_j) = 0$
for $0 \le s \le |H|-2$. That is,
$$
\tr A^s \prod_{j \not\in H} (A-\theta_j I) = 0\qquad (0 \le s \le |H|-2).
$$
It follows that $\prod_{j \not\in H} (A-\theta_j I)$,
when written on the basis $\{ A_i \mid 0 \le i \le d \}$,
does not contain $A_s$ for $0 \le s \le |H|-2$.
Moreover, the number of factors is $r=d+1-|H|$, so that $A_s$ does not occur either
when $s > r$. Therefore,  $\prod_{j \not\in H} (A-\theta_j I)$
is a linear combination of the $A_s$ with $|H|-1 \le s \le r$ and, hence, for some constants $\alpha_i$,
$$
\frac{n}{\tr (A_r A^{r})}\prod_{j \not\in H} (A-\theta_j I)=\sum_{s=|H|-1}^{r}\frac{\alpha_s}{k_s}A_s.
$$
Now comparing coefficients of $A^{r}$, we see that $\alpha_{r}=1$
(notice that $\tr (A_rA^{r})=nc_1\cdots c_{r}k_{r}$).
To obtain the value of $\alpha_i$ for $|H|-1\le i\le d-|H|$, multiply both terms of the above equation by $A_i$ and take traces. \qed

In the above we applied this twice, namely for
$d=2e+1$, $H=\{1,3,5,\ldots,d\}$, and for $d=2e$, $H=\{1,3,5,\ldots,d-1\}$.
In the latter case, \eqref{eqThm} yields the following expression for $z=-\alpha_{e-1}$.
$$
b_0b_1 \cdots b_e n z =
c_1c_2 \cdots c_{e+1} k_{e+1} n z =
- \tr (A_{e-1} \prod_{j \notin H} (A-\theta_j I)).
$$

Let us now take for $H$ the set of even indices, with or without 0.

\subsubsection*{\boldmath$H=\{0,2,4,\ldots,d\}$}
Let $d = 2e$ be even and suppose that $P_{id}$ takes the same value ($k_d$)
for all $i \in H = \{0,2,4,\ldots,d\}$. Then $|H| = e+1$, and
$\prod_{j \not\in H} (A-\theta_j I)$ is a multiple~of~$A_e$.
By Proposition \ref{0ddd} this happens if and only if
$\Gamma$ is antipodal with even diameter.
%
%

\subsubsection*{\boldmath$H=\{0,2,4,\ldots,d{-}1\}$}
Let $d = 2e+1$ be odd and suppose that $P_{id}$ takes the same value ($k_d$)
for all $i \in H = \{0,2,4,\ldots,2e\}$. Then $|H| = e+1$, and
$\prod_{j \not\in H} (A-\theta_j I)$ is a linear combination of $A_s$
for $e \le s \le e+1$.
By Proposition \ref{0ddd} this happens if and only if
$\Gamma$ is antipodal with odd diameter.

\subsubsection*{\boldmath$H=\{2,4,\ldots,d\}$}
Let $d = 2e$ be even and suppose that $P_{id}$ takes the same value
for all $i \in H = \{2,4,\ldots,d\}$. Then $|H| = e$, and
$\prod_{j \not\in H} (A-\theta_j I)$ is a linear combination of $A_s$
for $e-1 \le s \le e+1$.
As before we conclude that the $\theta_j$ with $j \not\in H$
are the eigenvalues of the array
$\{b_0,\ldots,b_{e-1};\,c_1,\ldots,c_{e-1},c_e{+}zb_e\}$
for some real $z \le 1$.
This time $d \in H$, and there is no conclusion about the sign of $z$.

For example, the Odd graph $O_5$ with intersection array
$\{5,4,4,3;\,1,1,2,2\}$ has eigenvalues $5$, $3$, $1$, $-2$, $-4$
and $P_{24} = P_{44}$. The eigenvalues 5, 3, $-2$ are those of
the array $\{5,4;\,1,-1\}$.

No primitive examples with $d > 4$ are known.
%
%
%

\subsubsection*{\boldmath$H=\{2,4,\ldots,d{-}1\}$}
Let $d = 2e+1$ be odd and suppose that $P_{id}$ takes the same value
for all $i \in H = \{2,4,\ldots,2e\}$. Then $|H| = e$, and
$\prod_{j \not\in H} (A-\theta_j I)$ is a linear combination of $A_s$
for $e-1 \le s \le e+2$.

For example, the Odd graph $O_6$ with intersection array
$\{6,5,5,4,4;\,1,1,2,2,\maysplit 3\}$ has eigenvalues
$6$, $4$, $2$, $-1$, $-3$, $-5$ and $P_{24} = P_{44}$.

No primitive examples with $d > 5$ are known.

\section*{Acknowledgments}
Part of this note was written while the second author was visiting
the Department of Combinatorics and Optimization (C\&O),
in the University of Waterloo (Ontario, Canada).
He sincerely acknowledges to the Department of C\&O the hospitality
and facilities received.
Also, special thanks are due to Chris Godsil for useful discussions
on the case of diameter four.

\end{document}